\documentclass[a4paper]{article}
\usepackage{latexsym}
\usepackage{amsmath,amssymb,amsfonts}

\newcommand{\noi}{\noindent}

\newtheorem{prethm}{{\bf Theorem}}

\newenvironment{thm}{\begin{prethm}{\hspace{-0.5
               em}{\bf.}}}{\end{prethm}}

\newtheorem{prepro}[prethm]{Proposition}
\newenvironment{pro}{\begin{prepro}{\hspace{-0.5
               em}{\bf.}}}{\end{prepro}}

\newtheorem{prelem}[prethm]{Lemma}
\newenvironment{lem}{\begin{prelem}{\hspace{-0.5
               em}{\bf.}}}{\end{prelem}}

\newtheorem{precor}[prethm]{Corollary}
\newenvironment{cor}{\begin{precor}{\hspace{-0.5
               em}{\bf.}}}{\end{precor}}

\newtheorem{preremark}{{\bf Remark}}

\newenvironment{rem}{\begin{preremark}\em{\hspace{-0.5
              em}{\bf.}}}{\end{preremark}}

\newtheorem{preexample}{{\bf Example}}

\newenvironment{example}{\begin{preexample}\em{\hspace{-0.5
               em}{\bf.}}}{\end{preexample}}

\newtheorem{preproof}{{\bf Proof.}}

\newenvironment{proof}[1]{\begin{preproof}{\rm
               #1}\hfill{$\Box$}}{\end{preproof}}

\newcommand{\rk}{{\rm rank}\,}

\date{}

\title{\bf\ Intersection matrices revisited}

\author{~ N. Ghareghani$^{\textrm{b}}$,~ E. Ghorbani$^{\textrm{c,b}}$, M. Mohammad-Noori$^{\textrm{a,}}$\thanks{Corresponding author}\\
{\small {\em $^{\rm a}$Department of Mathematics, Statistics and Computer Science, University of Tehran, Tehran, Iran}}\\
{\small {\em $^{\rm b}$School of Mathematics, Institute for Research in Fundamental Sciences (IPM),}}\\
 {\small {\em P.O. Box 19395-5746, Tehran, Iran}}\\
{\small {\em $^{\rm c}$Department of Mathematics,  K.N. Toosi University of Technology,  P. O. Box 16315-1618, Tehran, Iran}}
\\{\footnotesize Emails: {\tt
ghareghani@ipm.ir, e\_ghorbani@ipm.ir, morteza@ipm.ir,
mnoori@khayam.ut.ac.ir}} }
 \vspace{4mm}

\begin{document}
\maketitle

\begin{abstract}
\noi Several intersection matrices
 of $s$-subsets vs. $k$-subsets of a $v$-set are
introduced in the literature. We study these matrices
systematically through counting arguments and generating function
techniques. A number of new or known identities appear as natural
consequences of this viewpoint; especially, appearance of the
derivative operator $d/dz$ and some related operators reveals
some connections between intersection matrices and the
``combinatorics of creation-annihilation". As application, the
eigenvalues of several intersection matrices including some generalizations of
 the adjacency matrices of the Johnson scheme are derived; two new bases
for the Bose--Mesner algebra of the Johnson scheme are introduced
and the associated intersection numbers are obtained as well.
 Finally, we determine the rank of some intersection matrices.

\vspace{3mm}
\noindent {\em AMS Classification}:  05B20; 05E30\\
\noindent{\em Keywords}: Intersection matrices; Inclusion
matrices; derivative operator; Johnson scheme
 \end{abstract}

\section{Introduction}
Let $s$, $k$, and $v$ be integers satisfying $0\leq s, k\leq
v$. We fix a $v$-set $V$ throughout this paper. The inclusion
matrix $W_{s,k}$  is a $(0,1)$-matrix whose
rows and columns are indexed  by $s$-subsets and $k$-subsets of $V$, respectively,
and $W_{s,k}(S,K)=1$ if and only if
$S\subseteq K$. This matrix has interesting
properties and arises in many combinatorial problems, especially
in design theory and extremal set theory (see
\cite{Beth,fra,gkmm,w82,w84,w90}). It satisfies several nice
identities among which is
\begin{equation}\label{WWW} W_{i,s} W_{s,k}= {{k-i}\choose {s-i}}W_{i,k},\end{equation}
which holds for $i\le s\le k\le v$ (see \cite{KM, k,
w82,w84,w90}). Another $(0,1)$-matrix which is closely related to
$W_{s,k}$ is  the exclusion matrix $\overline{{W}}_{s,k}$  with the
same row and column indices as $W_{s,k}$ where the
$(S,K)$ is $1$ if and only if $S\cap K=\emptyset$.

Both of the inclusion and exclusion matrices may be regarded as
intersection matrices in the sense that the entry $(S,K)$ of them
only depends on $|S\cap K|$. Some other intersection matrices are
also studied in the literature in different contexts in
combinatorics, including design theory, association scheme and
extremal set theory,  the most significant properties of which
are the combinatorial identities they satisfy. The goal of the
present paper is to introduce and investigate a more general
framework in which several intersection matrices arise as
special cases and the identities involving them are derived more naturally. This uniform framework demonstrates
 the relation between intersection matrices and some operators of the form $\phi(z)\frac{d}{dz}$.
These operators were studied previously in \cite{Scherk} and more
recently in \cite{BF} and the references therein.

The paper is organized as follows. Section 2 contains basic properties of the derivative operator and some binomial identities which we will use later.
 In Section 3 we introduce several  intersection matrices and study relations between
 them.
 Particularly, we show that all these matrices can be
extracted from one, namely $F_{s,k,t}(z)$, a matrix with polynomial entries in variable $z$. We also show that studying
identities containing this matrix, produces identities containing
the other matrices.
In Section~4 we calculate the matrix
product $W^{\top}_{i,s} F_{i,k,t}$ as a linear combination of
derivatives of $F_{s,k,t}$. This reveals a close connection between
this matrix product and the operator $\frac{zd}{dz}$.
  Section 5 is the application section in which we
introduce two new bases for the Bose--Mesner algebra of the
Johnson scheme using the intersection matrices above.
  The eigenvalues of some generalizations of
 the adjacency matrices of the Johnson scheme are  also derived.
All the eigenvalues of these matrices can be expressed in terms of the polynomials $\mu_j(z)=\sum_{i=j}^{t}{k-j \choose i-j}{v-j-i\choose k-i}z^i$ for $j=0,1,\ldots,t$.


\section{Operators and basic notation}\label{opert}

Let $D$ denote the derivative operator $\frac{d}{dz}$, and
$(zD)_n$ denote the {\em falling factorial} $(zD)(zD-{\bf
1})\cdots(zD-(n-1){\bf 1})$, where ${\bf 1}$ is the identity
operator. For convenience we replace $\alpha {\bf 1}$ by
$\alpha$ like $(zD)(zD- 1)\cdots(zD-(n-1))$. Here are some of the identities containing the derivative operator:
\begin{enumerate}
  \item $(zD)^n =\sum_{k=1}^nS(n,k)z^kD^k$ for $S(n,k)$ the Stirling numbers of the second kind;
  \item $(zD)_n =z^nD^n$;
  \item $(zD-k)_n=n!\sum_{r=0}^n{n+k-r-1\choose n-r}(-1)^{n-r}\frac{z^r}{r!}D^r$;
  \item $D^\ell(z^rD^r)=\sum_{j=0}^\ell{\ell\choose
j}(r)_{\ell-j}z^{r-\ell+j}D^{r+j}$.
\end{enumerate}

  We frequently make use of the following binomial identities (see Chapter 5 of \cite{knuth}):
\begin{align}
\label{knuth0} & \sum_k (-1)^k{\ell \choose m+k}{s+k \choose
n}=(-1)^{\ell+m}{s-m \choose n-\ell},~~~\ell\geq 0,~~\hbox{and}\\
\label{knuth} & \sum_{k\leq \ell} (-1)^{k} {\ell-k\choose m}{s
\choose k-n}= (-1)^{\ell+ m}{s-m-1\choose
\ell-m-n},~~~\ell,m,n\geq 0.
\end{align}

The coefficient of $z^i$ in a polynomial (or a generating
function) $p(z)$ is denoted as $[z^i]p(z)$. For two matrices $A$
and $B$, we write $A \equiv B$ if $A$ can be obtained from $B$ by
a permutation of the rows and a permutation of the columns. For
instance $\overline{W}_{sk}\equiv W_{s,v-k}$; this is because for
$S\in{V\choose s}$ and $K\in{V\choose k}$,  $|S\cap K|=\emptyset$
if and only if $S\subseteq V\setminus K$.

\section{Intersection matrices}
 Let $s, k$ and $v$ be
integers satisfying $0\leq s\leq k\leq v$. Let ${V\choose s}$ and
${V\choose k}$ denote the sets of $s$-subsets and $k$-subsets of
$V$, canonically ordered somehow. Then an {\it intersection
matrix} (relative to this setup) is a ${v\choose
s}\times{v\choose k}$ matrix with $(S,K)$ entry as a function of
$|S\cap K|$ (but not otherwise dependent on the specific subsets
$S$ and $K$). The inclusion matrix $W$ is one such since
$${|S\cap K|\choose s}=\left\{\begin{array}{ll} 1 & \hbox{if $S\subseteq K$,} \\ 0 & \hbox{otherwise.}\end{array}\right.$$
The exclusion matrix $\overline{W}$ is another such since
$${s-|S\cap K|\choose s}=\left\{\begin{array}{ll} 1 & \hbox{if $S\cap K=\emptyset$,} \\ 0 & \hbox{otherwise.}\end{array}\right.$$

\noindent{\bf Definitions.} The intersection matrices considered
herein are:
\begin{enumerate}
  \item $A_{s,k,t}$ with $(S,K)$ entry ${|S\cap K|\choose t}$,
   a generalization of $W_{s,k}=A_{s,k,s}$;
  \item $F_{s,k,t}(z)$ with $(S,K)$ entry $\sum_{i=0}^t{|S\cap K|\choose i}z^i$;
  \item $U_{s,k,t,\ell}$ with $(S,K)$ entry $\sum_{i=0}^t(-1)^{i-\ell}{i\choose\ell}{|S\cap K|\choose i}={|S\cap K|\choose\ell}{t-|S\cap K|\choose t-\ell}$;
  \item $N_{s,k,t}$ with $(S,K)$ entry $\sum_{i=0}^t(-1)^{t+i}{|S\cap K|\choose i}={|S\cap K|-1\choose t}$;
  \item $F_{s,k}(z)=F_{s,k,s}(z)$;
  \item $U_{s,k,\ell}= U_{s,k,s,\ell}$
   which is a $(0,1)$-matrix whose $(S,K)$ entry is 1 if and only if  $|S\cap K|=\ell$.
\end{enumerate}

The obvious relations among these include
\begin{enumerate}
  \item $F_{s,k,t}(z)=\sum_{i=0}^t A_{s,k,i}z^i$;
  \item $F_{s,k,t}(z)=\sum_{\ell=0}^t U_{s,k,t,\ell}(z+1)^\ell$;
  \item $A_{s,k,i}=\frac{D^i}{i!}F_{s,k,t}(z)|_{z=0}$;
  \item $U_{s,k,t,\ell}=\frac{D^\ell}{\ell!}F_{s,k,t}(z)|_{z=-1}$;
  \item $F_{s,k}(z)=(z+1)^{|S\cap K|}$;
  \item $U_{s,k,t,\ell}=\sum_{i=\ell}^t(-1)^{i-\ell}{i\choose\ell}A_{s,k,i}$;
  \item $A_{s,k,i}=\sum_{\ell=i}^t{\ell\choose i}U_{s,k,t,\ell}$;
  \item $N_{s,k,t}=(-1)^tU_{s,k,t,0}$;
  \item $N_{s,k,t}=\sum_{i=0}^t(-1)^{t-i}A_{s,k,i}$;
  \item $A_{s,k,t}=N_{s,k,t}+N_{s,k,t-1}$.
\end{enumerate}

The following properties of the intersection matrices are
straightforward.
\begin{pro}

\label{propert}
\begin{itemize}
\item[\rm(i)]  $A_{t,v-k,t}\equiv
\overline{W}_{t,k}$ and $A_{v-t,k,k}\equiv \overline{W}_{t,k}$.
\item[\rm(ii)] $U_{t,k,t}=W_{t,k}$ and $U_{k,t,t}=W^{\top}_{t,k}$.
\item[\rm(iii)] If $t\geq \min(s,k)$, then $F_{s,k,t}=F_{s,k}$ and $U_{s,k,t,\ell}=U_{s,k,\ell}$.
\item[\rm(iv)] $U_{t,k,0}=\overline{W}_{t,k}$.
\item[\rm(v)] $A_{s,k,t}={W}_{t,s}^\top W_{t,k}$.
\item[\rm(vi)] $F_{s,k,t}^\top=F_{k,s,t}$.
Hence, $F_{k,k,t}$ is a symmetric matrix.
\item[\rm(vii)] $U_{s,k,t,\ell}(S,K)\neq 0$ only if $|S\cap K|\in B$ where $B=\{\ell\}\cup\{t+1,t+2,\ldots,\min(s,k)\}$.
\item[\rm(viii)] If $\ell\le t\le\min(s,k)$, then the number of nonzero elements in each row of $U_{s,k,t,\ell}$ is $\sum_{i\in B} {s \choose i} {v-s\choose k-i}$.
\item[\rm(ix)]$U_{s,k,t,\ell}= \sum_{i\in B} {i\choose \ell}{t-i\choose t-\ell}  U_{s,k,i}$.
\item[\rm(x)] There are exactly ${{v-s}\choose k}+ \sum_{i=t+1}^{\min(s,k)}{s \choose i} {v-s\choose k-i}$ nonzero elements in each row of $N_{s,k,t}$.
\item[\rm(xi)] There are exactly ${{v-s}\choose k}+{s \choose k}$ nonzero elements in each row of $N_{s,k,k-1}$.
\item[\rm(xii)] There are exactly ${{v-s}\choose k}+{v-s \choose k-s}$ nonzero elements in each row of $N_{s,k,s-1}$.
\end{itemize}
\end{pro}

\begin{rem}
\label{N13-14}
The matrix $N_{s,k,t}$ was introduced in \cite{KMT} and discussed further in \cite{TEZ} as an
auxiliary tool to speed up an algorithmic search for finding $t$-designs.
 More precisely, the matrices $N_{7,7,6}$ with $v=14$ and $N_{6,6,5}$ with $v=13$ had important roles in finding $6$-$(14,7,4)$ designs (see \cite{KMT}). It was observed that the small number of nonzero elements in each row of $N_{s,k,t}$ (which is obtained by choosing proper values of $s$ as it can be seen from Proposition~\ref{propert}) is a useful property for this. The matrix $N_{7,7,6}$ with $v=14$ has only $2$ nonzero elements in each of its rows and the matrix $N_{6,6,5}$ with $v=13$ has $8$ nonzero elements in each of its rows.
\end{rem}

\begin{pro}\label{NewFWAU} The followings hold:
\begin{itemize}
\item[\rm(i)] $U_{s,k,\ell} \equiv U_{s,v-k,s-\ell}\equiv U_{v-s,k,k-\ell} \equiv U_{v-s,v-k,v-s-k+\ell}$;
\item[\rm(ii)] $F_{v-s,k}(z) \equiv (z+1)^k F_{s,k}(\frac{-z}{z+1})$;
\item[\rm(iii)] $F_{s,v-k}(z) \equiv (z+1)^s F_{s,k}(\frac{-z}{z+1})$;
\item[\rm(iv)] $F_{v-s,v-k}(z) \equiv (z+1)^{v-s-k}F_{s,k}(z)$.
\end{itemize}
\end{pro}
\begin{proof}
{Note that $|S\cap K|=\ell$ if and only if $|S\cap (V\setminus
K)|=s-\ell$. Thus, we have $U_{s,k,\ell} \equiv
U_{s,v-k,s-\ell}$. The rest of part (i) is proved similarly.
Since the proofs of the remaining parts are similar, we only prove
(ii):
\begin{align*}
F_{v-s,k}(z)(V\setminus S,K)&= (z+1)^{|(V\setminus S)\cap K|}\\
&= (z+1)^{k-|S\cap K|}\\
 &=(z+1)^{k} (1-\tfrac{z}{z+1})^{|S\cap K|}\\
 &=(z+1)^{k} F_{s,k}(\tfrac{-z}{z+1})(S,K).
 \end{align*}}
\end{proof}

\section{Matrix products}

The following theorem gives some useful identities, among
them (\ref{UUabc}) and (\ref{AAabc}) are new. Moreover, two
new proofs are given for the known identity (\ref{ww'}) appeared first in \cite{w82}.

\begin{thm}\label{FWAU} The followings hold:
\begin{eqnarray}
U_{a,b,i} U_{b,c,j} &=& \sum_{\ell=0}^{\min(a,c)} \sum_{n=0}^{\ell}{\ell \choose n}{c-\ell \choose j-n}{a-\ell \choose i-n}{v-a-c+\ell \choose b-i-j+n}U_{a,c,\ell},\label{UUabc}\\
\label{WWt} W_{a,k}W_{b,k}^\top &=& \sum_{n=0}^{\min(a,b)}{v-b-a\choose v-k-n}A_{a,b,n} \label{ww'},\\
  A_{a,b,i} A_{b,c,j} &=& \sum_{n=0}^{\min(i,j)}{a-n \choose i-n} {c-n \choose j-n} {v-i-j \choose b+n-i-j}A_{a,c,n} \label{AAabc}.
\end{eqnarray}
\end{thm}

\begin{proof}{
For any given $a$-subset $A$ and  $c$-subset $C$, using simple
counting arguments, the entry $(A,C)$ of the matrix product
$U_{a,b,i} U_{b,c,j}$ is calculated as
$$\left( U_{a,b,i} U_{b,c,j}\right) (A,C)=\left|\left\{B\subseteq \{1,\ldots,v\} :|B|=b,~|B\cap A|=i,~ |B\cap C|=j \right\}\right|.$$
  To count the number of $b$-sets $B$ with the above constraints, let $\ell= |A\cap C|$ and $n= |A \cap B \cap C|$. To construct $B$, one should select $n$ points from $A\cap C$,
  $i-n$ points from $A\setminus C$, $j-n$ points from $C\setminus A$ and $b-i-j+n$ points from $A'\cap C'$. Hence,
  \begin{equation*}
  (U_{a,b,i} U_{b,c,j}) (A,C)=\sum_{\ell=0}^{\min(a,c)}\delta_{|A\cap C|,\ell}\sum_{n=0}^{\ell}
{\ell \choose n}{c-\ell \choose j-n}{a-\ell \choose
i-n}{v-a-c+\ell \choose b-i-j+n},
 \end{equation*}
  which proves (\ref{UUabc}).

To prove (\ref{ww'}), note that $W_{a,k} \equiv W^\top_{v-k,v-a}$
because $A \subseteq K$ if and only if $ V\setminus K \subseteq
V\setminus A$. If we apply simultaneously the same permutation on
the columns of $W_{a,k}$
 and on the rows of $W^\top_{b,k}$, then we see that
$W_{a,k}W^\top_{b,k} \equiv W^\top_{v-k,v-a} W_{v-k,v-b}$. Hence
\begin{align*}
W_{a,k}W^\top_{b,k} &\equiv  A_{v-a,v-b,v-k}\\
 &=[z^{v-k}]F_{v-a,v-b} \\
  &\equiv[z^{v-k}]\big( (z+1)^{v-a-b}F_{a,b} \big)\\
 &=\sum_{n=0}^{\min(a,b)} {{v-a-b}\choose {v-k-n}} [z^n]F_{a,b}\\
 &=\sum_{n=0}^{\min(a,b)} {{v-a-b}\choose {v-k-n}} A_{a,b,n},
 \end{align*}
concluding (\ref{ww'}). (We notice that the same ordering is used
for the rows and the columns of both matrices.) An alternative way to prove
(\ref{ww'}), is as follows:
\begin{align*}
W_{a,k}W^\top_{b,k}&=U_{a,k,a}U_{k,b,b}\\
&=\sum_{\ell}{v-a-b+\ell \choose k-a-b+\ell}U_{a,b,\ell}~~~~\hbox{(by (\ref{UUabc}))}\\
&=\sum_{\ell}{v-a-b+\ell \choose v-k}\sum_{n} (-1)^{n-\ell} {n \choose \ell} A_{a,b,n}\\
&=\sum_{n}(-1)^n A_{a,b,n}  \sum_{\ell}(-1)^{\ell}{n \choose
\ell}{v-a-b+\ell \choose v-k}\\
&=\sum_{n}{v-a-b \choose v-k-n}A_{a,b,n}~ ~\,\, ~~~~\hbox{(by
(\ref{knuth0}))}.
\end{align*}

 Now, replacing $a$ by $i$ and $k$ by $j$ in
(\ref{ww'}) and multiplying the identity from left and right
 by $W_{i,a}^\top$ and $W_{j,c}$, respectively, and using (\ref{WWW}), we conclude  (\ref{AAabc}).}
\end{proof}

We now turn to calculate the matrix product
$W_{s,j}^{\top}F_{j,k,t}$. This calculation reveals the
relationship between this matrix product and the derivative operator.
 The general form of these identities contains expressions of the
form $(\phi(z)D+\rho(z))^n$; such expressions are studied firstly
by Scherk in \cite{Scherk} and extensively by some authors in
recent years (see \cite{BlskThz} and the references therein).
Recently, many interesting properties associated with the algebra
of two operators, $A$ and $B$, satisfying $AB-BA={\bf 1}$  in
quantum physics. This relation is called  the {\it
creation-annihilation} axiom. A simple representation of this
algebra is obtained by taking $A=D$ and $B=z$ (see Section 2.4 of
\cite{BlskThz}). Also in \cite{BF} systematic evaluation of
expressions of the form $(\phi(z)D+\rho(z))^n$ using several
combinatorial models involving set partitions, permutations,
increasing trees and weighted lattice paths is studied. This is
discussed more in Remark \ref{fd} at the end of this section.

\begin{thm} The following identities hold:
\label{recIdnt}
\begin{enumerate}
  \item[\rm(i)] $W_{s-1,s}^\top F_{s-1,k,t}(z)= s F_{s,k,t}(z)-z DF_{s,k,t}(z)$,
 \item[\rm(ii)] $W_{s-1,s}^\top U_{s-1,k,t, \ell}= (s-\ell) U_{s,k,t,\ell}+(\ell+1) U_{s,k,t,\ell+1}$.
\end{enumerate}
\end{thm}
\begin{proof}{ By Proposition~\ref{propert}(i) and (\ref{WWW}), we have
$W_{s-1,s}^\top A_{s-1,k,i}=(s-i)A_{s,k,i}$. Therefore
\begin{align*} W_{s-1,s}^\top F_{s-1,k,i}(z)&= \sum_{i=0}^t (s-i)A_{s,k,i} z^i\\
&= s \sum_{i=0}^t A_{s,k,i} z^i-z \sum_{i=1}^ti A_{s,k,i} z^{i-1},
 \end{align*}
which proves (i). By applying the operator $D^{\ell}$ to (i) we
get
$$ W_{s-1,s}^\top D^{\ell}F_{s-1,k,i}(z)= (s-\ell) D^{\ell} F_{s,k,i}(z)-z D^{\ell+1} F_{s,k,i}(z). $$
Setting $z=-1$ we obtain (ii). }\end{proof}

\begin{rem}\label{cncludEqu} From Theorem~\ref{recIdnt}, some new identities can be derived.
By setting $t=s$ in (i) and using Proposition~\ref{propert}(ii),
\begin{equation*}
W_{s-1,s}^\top F_{s-1,k}(z)= s F_{s,k}(z)-z DF_{s,k}(z).
\end{equation*}
Transposing the two sides of this equation and exchanging $k$ and $s$, we
have
\begin{align*}
F_{s,k-1,t}(z)W_{k-1,k}&= k F_{s,k,t}(z)-z DF_{s,k,t}(z),~~\hbox{and} \\
F_{s,k-1}(z)W_{k-1,k}&= k F_{s,k}(z)-zDF_{s,k}(z).
\end{align*}
\end{rem}

In  Theorem~\ref{recIdnt}(ii), the expression $W_{i,s}^\top
U_{i,k,t,\ell}$ is calculated for $i=s-1$, but how can we
calculate this expression in general? In the special case of
$t=k$,
 the answer is simply obtained by using (\ref{UUabc}) as follows:
 $$W_{i,s}^\top U_{i,k,\ell}= U_{s,i,i} U_{i,k,\ell} =\sum_{h=\ell}^s {h
\choose \ell} {s-h \choose i-\ell} U_{s,k,h}.$$
 The following theorem gives the answer in general.

\begin{thm} \label{WFUU} Let $L_{s,i}(z)=\sum_{r=0}^{s-i}(-1)^r{s-r\choose i}\frac{z^r}{r!}D^r.$ Then the following identities hold:

\begin{itemize}
\item[\rm(i)]${\displaystyle W_{i,s}^\top
F_{i,k,t}(z)=L_{s,i}(z)F_{s,k,t}(z),}$
 \item[\rm(ii)]${\displaystyle W_{i,s}^\top U_{i,k,t,\ell}=\sum_{h=\ell}^{\ell+s-i}{h\choose
\ell}{s-h\choose i-\ell}U_{s,k,t,h}.}$
\end{itemize}
\end{thm}
\begin{proof}{
(i) By Theorem~\ref{recIdnt}(i), $W_{s-1,s}^\top F_{s-1,k,t}(z)=
(s-zD) F_{s,k,t}(z)$. On the other hand, from
(\ref{WWW}) it follows that
\begin{equation*}
   W_{i,i+1}W_{i+1,i+2}\cdots W_{s-1,s}=(s-i)!W_{i,s}.
\end{equation*}
Now by iterative use of Theorem~\ref{recIdnt}(i), we have
\begin{align*}
    W_{i,s}^{\top} F_{i,k,t}&=\frac{1}{(s-i)!}W^{\top}_{s-1,s} W^{\top}_{s-2,s-1}\cdots
    W^{\top}_{i,i+1} F_{i,k,t}\\
    &=\frac{1}{(s-i)!}(s-zD)(s-1-zD)\cdots (i+1-zD)F_{s,k,t}.
\end{align*}
By the property 3 of Section~2, the operator
$\frac{1}{(s-i)!}(s-zD)(s-1-zD)\cdots (i+1-zD)$ in the last
expression can be simplified as
$L_{s,i}(z)=\sum_{r=0}^{s-i}(-1)^r{s-r\choose
i}\frac{z^r}{r!}D^r.$

(ii) Applying the operator $D^\ell/\ell!$ on (i), we have
\begin{align*}
 W_{i,s}^\top \frac{D^\ell}{\ell!}
F_{i,k,t}&=\frac{1}{\ell!}\sum_{r=0}^{s-i}\frac{(-1)^r}{r!}{s-r\choose
i}D^\ell(z^rD^r) F_{s,k,t}\\
&=\frac{1}{\ell!}\sum_{r=0}^{s-i}\frac{(-1)^r}{r!}{s-r\choose
i}\sum_{j=0}^\ell{\ell\choose
j}(r)_{\ell-j}z^{r-\ell+j}D^{r+j}F_{s,k,t}(z).
\end{align*}


 Letting $z=-1$, we get
\begin{align*}
   W_{i,s}^\top U_{i,k,t,\ell} &=\sum_{r=0}^{s-i}\sum_{j=0}^\ell{s-r\choose
i}{r+j\choose j}{r\choose\ell-j}(-1)^{\ell-j}U_{s,k,t,r+j}\\
&=\sum_{h=0}^{\ell+s-i}a_h U_{s,k,t,h},
\end{align*}
where $a_h=\sum_{j=0}^h(-1)^{\ell-j}{h\choose
j}{h-j\choose\ell-j}{s-h+j\choose i}$. From (\ref{knuth}) it
follows that $a_h={h\choose \ell}{s-h\choose i-\ell}$.}
\end{proof}

\begin{rem}\label{fd} The calculations in this section show some connections between intersection
matrices and operators of the form $(\phi(z)D)^n$ (or more
generally $(\phi(z)D+\rho(z))^n$). These operators are studied in
\cite{BlskThz} and in more details in \cite[Section 6.3]{BF}.
There are more such connections. For instance, we can prove that
\begin{equation} \label{wjf0} W_{s,j}F_{j,k}(z)= (z+1)^{-v+j+k}\sum_{r=0}^{j-s}(-1)^r{v-s-r \choose v-j}\frac{z^rD^r}{r!}((z+1)^{v-s-k}F_{s,k}(z)).
\end{equation}
On the other hand,  by (\ref{WWt}),
$W_{j-1,j}F_{j,k}(z)=(\phi(z)D+\rho_j(z))F_{j-1,k}(z)$ where
$\phi(z)=-z^2-z$ and $\rho_j(z)=kz+v-j+1$. Using techniques
similar to the ones used in the proof of Theorem~\ref{WFUU}, it
follows that
\begin{equation*}
W_{s,j}F_{j,k}(z)=\frac{1}{(j-s)!}(\phi(z)+\rho_{s+1})(\phi(z)+\rho_{s+1}-1)\cdots(\phi(z)+\rho_{s+1}-j+s+1)
F_{s,k}(z).
\end{equation*}
We believe that deeper relations of this form help one in
studying more useful properties of intersection matrices. Such
operators were studied firstly in \cite{Scherk} and more recently
in \cite{BF}.
 \end{rem}

\section{Some applications}
This section contains some applications of  the results obtained
so far. One important application is
deriving the eigenvalues of the matrices $F_{k,k,t}(z)$, $U_{k,k,\ell,t}$ and $N_{k,k,t}$ based on Wilson's
method for computing the eigenvalues of $A_{k,k,i}$ \cite{w82} (cf. \cite{w84,w84b}).
 To the
best of our knowledge, the rest of the results of this section
are new. Among which are introducing two new bases for the
Bose--Mesner algebra of the Johnson scheme and obtaining the
associated intersection numbers.
 We also determine the rank of some intersection matrices.

\subsection{Johnson scheme}

An {\em association scheme with $d$ classes} is a set of $d+1$
square $(0,1)$-matrices $X_0,X_1,\ldots,X_d$ which satisfy
\begin{enumerate}
  \item[\rm(i)] $\sum_{i=0}^dX_i=J$,
  \item[\rm(ii)] $X_0=I$,
  \item[\rm(iii)] $X_i=X_i^\top$, for $i=0,1,\ldots,d$,
  \item[\rm(iv)] $X_iX_j=\sum_{\ell=0}^d a_{ij}^\ell X_\ell$, for $i,j\in\{0,1,\ldots,d\}$.
\end{enumerate}
The numbers $a_{ij}^\ell$ are called the {\em intersection
numbers} of the association scheme. From (i) we see that the
matrices $X_i$ are linearly independent, and by use of (ii)--(iv)
we see that they generate a commutative $(d+1)$-dimensional
algebra of symmetric matrices with constant diagonal. This algebra
is called the {\em Bose--Mesner algebra} of the association scheme.

A Bose--Mesner algebra has a basis $\{E_0 = \frac{1}{n} J ,
E_1,\ldots, E_d \}$ of idempotents, that is,
 $E_i E_j = \delta_{i, j} E_i$ where $\delta_{i, j}$ is the Kronecker symbol.
The change-of-coordinates matrix $P = [p_{i j}]$ defined by $X_j
= \sum_i p_{i j} E_i$ has the property that
 $p_{i j}$ is an eigenvalue of $X_j$ whose eigenspace is the column space of $E_i$.
The matrix $P$ of eigenvalues  contains many properties of the
scheme from which many parameters of the scheme (such as
$a_{ij}^\ell$, etc.) can be obtained (see \cite{d}). In this
regard, the eigenvalues of different bases of an association
scheme are important subjects and worth to study.

The {\em Johnson scheme} $J(v,k)$ is a $k$-class association
scheme in which the rows and the columns of each $X_i$ is indexed
by all  $k$-subsets of a $v$-set and $X_i(K_1,K_2)=1$ if and only
if $|K_1\cap K_2|=k-i$, for $i=0,1,\ldots,k$. In other words,
$X_i=U_{k,k,k-i}$. In this section we introduce two new bases for
the Bose--Mesner algebra of $J(v,k)$ and obtain the associated
intersection numbers.

The first new basis for the Bose--Mesner algebra of $J(v,k)$ is
$\left\{A_{k,k,i} : i=0,\ldots,k\right\}$; this follows from the
identities $U_{s,k,t,\ell}=\sum_{i=\ell}^t (-1)^{i-\ell}
{i\choose \ell} A_{s,k,i}$ and $A_{s,k,i}=\sum_{\ell=i}^t  {\ell
\choose i} U_{s,k,t,\ell}$.

To introduce the second basis we define the matrix
$B_{s,k,\ell}$  as
$$B_{s,k,\ell}(S,K)=\left\{\begin{array}{ll}1 & \hbox{if $|S\cap K|\ge\ell$,} \\0 & \hbox{otherwise.}\end{array}\right.$$
Whence, we have $B_{s,k,\ell}=\sum_{\ell'=\ell}^{s}U_{s,k,\ell'}$
and $U_{s,k,\ell}=B_{s,k,\ell}-B_{s,k,\ell+1}$. This shows that
the matrices $\left\{B_{k,k,\ell}: \ell=0,\ldots,k\right\}$ form
a basis for the Bose--Mesner algebra of $J(v,k)$. The relation
between the two new bases is demonstrated  below.
\begin{pro}\label{UgeqL}
If $\ell>0$, then
\begin{equation}\label{U>}
  B_{s,k,\ell}= \sum_{i=\ell}^{s} (-1)^{i-\ell}{{i-1} \choose {\ell-1}}A_{s,k,i}.
\end{equation}
\end{pro}
\begin{proof}{
Let $G_{s,k}(z)=F_{s,k}(z-1)$,
$H_{s,k}(z)=\frac{1}{z-1}(G_{s,k}(z)-G_{s,k}(1))$ and
$G^{+}_{s,k}(z)=\sum_{\ell} B_{s,k,\ell}z^{\ell}$. Then
$G_{s,k}(z)=\sum_{\ell=0}^s U_{s,k,\ell} z^{\ell}$ and
$H_{s,k}(z)=\sum_{i=1}^s A_{s,k,i}(z-1)^{i-1}$. Moreover, by
$U_{s,k,\ell}=B_{s,k,\ell}-B_{s,k,\ell+1}$, we have
\begin{align*}
(z-1)G^{+}_{s,k}(z)&=zG_{s,k}(z)-G_{s,k}(1)\\
&=z\left(G_{s,k}(z)-G_{s,k}(1)\right)+(z-1)G_{s,k}(1)\\
&=z(z-1)H_{sk(z)}+(z-1)G_{s,k}(1).
\end{align*}

Hence, $G^{+}_{s,k}(z)=zH_{s,k}(z)+G_{s,k}(1)$ and for $\ell>0$ we
get
\begin{align*}
B_{s,k,\ell}&=[z^{\ell}]G^{+}_{s,k}(z)\\
&=[z^{\ell-1}]H_{s,k}(z)\\
&=D^{\ell-1}\left(\sum_{i=1}^s
A_{s,k,i}(z-1)^{i-1}\right)\big|_{z=0}\\
&=\sum_{i=\ell}^{s} (-1)^{i-\ell}{{i-1} \choose
{\ell-1}}A_{s,k,i}.
\end{align*}
}\end{proof}


Define the intersection numbers $r_{ij}^\ell$ and $p_{ij}^\ell$ as
$$A_{k,k,i}A_{k,k,j}=\sum_{\ell=0}^k r_{ij}^\ell
A_{k,k,\ell}, ~\hbox{and}~~ U_{k,k,i}U_{k,k,j}=\sum_{\ell=0}^k
p_{ij}^\ell U_{k,k,\ell}.$$
 From (\ref{AAabc}) and (\ref{UUabc}) it follows that:
\begin{pro} \label{Scmp} The values of intersection numbers $r_{ij}^\ell$ and $p_{ij}^\ell$ are as follows:
$$r_{ij}^\ell={v-i-j\choose k-i-j+\ell}
  {k-\ell\choose i-\ell}{k-\ell\choose j-\ell},
  ~\hbox{and}~~
   p_{ij}^\ell=\sum_{e=0}^\ell
{\ell \choose e} {k-\ell \choose i- e} {k-\ell \choose
j-e}{v-2k+\ell \choose k-i-j+e}.$$
\end{pro}

\subsection{Eigenvalues and rank of intersection matrices}

 The eigenvalues of $U_{k,k,k-\ell}$ (for $\ell=0,1,\ldots,k$), the adjacency matrices of the Johnson scheme $J(v,k)$, can be
expressed in terms of ``Eberlein polynomials'' (see \cite{bi,d})
which are
$$ E_j=\sum_{i=0}^\ell(-1)^{\ell-i}{k-i\choose \ell-i}{k-j\choose
i}{v-k+i-j\choose i},$$ with multiplicity ${v\choose j}-{v\choose
j-1},$  for $j=0,1,\ldots,k.$ In this section, we
 obtain the
eigenvalues of $F_{k,k,t}(z)$ and $U_{k,k,t, \ell}$ as well as
$B_{k,k,\ell}$. The eigenvalues of $W^\top_{s,k}F_{s,k,t}(z)$ and
$W_{s,k}^\top U_{s,k,\ell}$ are also determined. Moreover, we
give a closed form for the eigenvalues and the rank of
$N_{k,k,k-1}$. The rank of $U_{t,k,\ell}$ is also investigated.
It is interesting that all the eigenvalues of above matrices can be expressed in terms of the polynomials $
   \mu_j(z)=\sum_{i=j}^{t}{k-j \choose i-j}{v-j-i\choose k-i}z^i$ for $j=0,1,\ldots,t$.

The following lemma which gives the eigenvalues and the
corresponding eigenvectors of $A_{k,k,i}$ was proved by Wilson
\cite{w82}
 with a proof based on Equation (\ref{AAabc}). The following
decomposition of $\mathbb{R}^{v \choose k}$ is used in \cite{w82}:
fix $k$ and let $R_j$ denote the row-space of $W_{j,k}$ over the
field $\mathbb{R}$. From (\ref{WWW}) it follows that $R_0\subseteq R_1\subseteq
\cdots \subseteq R_k=\mathbb{R}^{v \choose k}.$ Now let
$V_0=R_0$, and $V_j:=R_j\cap R_{j-1}^\perp$ for $j=1,\ldots,k$.
Then $\mathbb{R}^{v\choose k}=V_0\oplus V_1\oplus \cdots\oplus
V_k$, and $V_j$ has dimension ${v\choose j}-{v\choose j-1}$. We note that, as it is well known, if $s\leq k\leq v-s$, then $\rk W_{sk}={v \choose s}$ (see \cite{de caen, gll, gju, w90}); moreover, an explicit right inverse for $W_{sk}$ in this case is given in \cite{gll,K,bap}.

\begin{lem}\label{idA}
{\rm (\cite{w82})}
With the above definitions, for any ${\bf x}\in V_j$, $A_{k,k,i}
{\bf x}^{\top}=\lambda_j {\bf x}^{\top}$, where
$$\lambda_j=\left\{
  \begin{array}{ll}{v-i-j\choose k-i}{k-j\choose i-j} & \hbox{if $i\geq j$,} \\0 & \hbox{otherwise.}\end{array}\right.$$
   In other words, the vectors of $R_i^{\perp}$ are
eigenvectors corresponding the eigenvalue $0$ and the vectors in
$V_j$, for $j=0,\ldots,i$ are eigenvectors corresponding the
eigenvalue ${k-j\choose i-j}{v-j-i\choose k-i}$.
\end{lem}

The following theorem determines the eigenvalues of
$F_{k,k,t}(z)$. Before that we need  further definitions: fix $k$
and let $R_j(z)$ denote the row-space of $W_{j,k}$ over the field
of rational functions $\mathbb{R}(z)$ and let $V_0(z):=R_0(z)$,
$V_j(z):=R_j(z)\cap R_{j-1}(z)^\perp$ for $j=1,\ldots,k$. Note
that a basis of $R_j$ (resp. $V_j$) over the ground field
$\mathbb{R}$ is also  a basis for $R_j(z)$ (resp. $V_j(z)$) over
the field $\mathbb{R}(z)$.

\begin{thm}\label{teigF} Let $0\leq t\leq k\leq v/2$. Consider $F_{k,k,t}(z)$ as a matrix with entries in the
field of rational functions $\mathbb{R}(z)$. Then the eigenvalues
of $F_{k,k,t}(z)$ are
$$\mu_0(z)^{{v \choose 0}},\mu_1(z)^{{v \choose 1}-{v \choose 0}},\ldots,\mu_{t}(z)^{{v \choose {t}}-{v \choose {t-1}}},0^{{v \choose k}-{v \choose {t}}},$$
where the exponents indicate the multiplicity and
\begin{equation}\label{eigF}
   \mu_j(z)=\sum_{i=j}^{t}{k-j \choose i-j}{v-j-i\choose
k-i}z^i,
\end{equation}
for $j=0,1,\ldots,t$. Furthermore, with the above notations, the
vectors in $V_j(z)$ are eigenvectors corresponding to $\mu_j$,
for $j=0,\ldots,t$. The vectors of $R_t(z)^{\perp}$ are
eigenvectors corresponding to the eigenvalue $0$.
\end{thm}

\begin{proof}{Considering $F_{k,k,t}(z)=\sum_{i=0}^t A_{k,k,i} z^i$, the proof follows from Lemma~\ref{idA}.}
\end{proof}
Now, it is easily seen that Eberlin polynomials, defined at the
beginning of this section, are obtained from polynomials
$\mu_j(z)$, in the case $t=k$ as follows:
\begin{equation*}
E_j=\frac{D^{k-\ell}}{(k-\ell)!}\mu_j(z)\big|_{z=-1}.
\end{equation*}
In general we have
\begin{cor}\label{teigU} Let $0\leq t\leq k\leq v/2$. The eigenvalues of $U_{k,k,t, \ell}$ are
$$\lambda_0^{{v \choose 0}},\lambda_1^{{v \choose 1}-{v \choose 0}},\ldots,\lambda_{t}^{{v \choose {t}}-{v \choose {t-1}}},0^{{v \choose k}-{v \choose {t}}},$$
where
\begin{equation}\label{eigU}
   \lambda_j=\frac{D^{\ell}}{\ell!}\mu_j(z)\big|_{z=-1}=\sum_{i=\ell}^{t}(-1)^{\ell+i}{i\choose \ell}{k-j\choose i-j}{v-j-i\choose k-i},
\end{equation}
for $j=0,1,\ldots,t$.
\end{cor}

In the previous subsection we saw that $\{B_{k,k,\ell}: \, \ell=0,
\ldots, k\}$ gives a new basis for the Johnson scheme, so it
is important to calculate their eigenvalues which are given in
the following corollary.

\begin{cor} Let $\ell\geq 0$ and $k \leq v/2$ and let $\mu_j(z)$ be as in (\ref{eigF}) with the additional condition $t=k$. Moreover, let
$$\nu_j(z)=\frac{z\mu_j(z-1)-\mu_j(0)}{z-1},~~~~j=0,\cdots,k.$$
 The eigenvalues of $B_{k,k,\ell}$ are $\frac{D^{\ell}}{\ell!}\nu_j(z)|_{z=0}$
with multiplicity ${v\choose j}-{v\choose j-1}$, for
$j=0,1,\ldots,k$.
 Furthermore, if $\ell>0$, then 
$$
   \frac{D^{\ell}}{\ell!}\nu_j(z)\big|_{z=0}=\sum_{i=\ell}^k(-1)^{\ell+i}{i-1\choose \ell-1}{k-j\choose i-j}{v-j-i\choose
k-i}.
$$
\end{cor}
\begin{proof}
{Considering the notation used in the proof of
Proposition~\ref{UgeqL}, we have
\begin{equation*}
G_{k,k}^+(z)=\frac{z}{z-1}(F_{k,k}(z-1)-F_{k,k}(0))+F_{k,k}(0).
\end{equation*}
Now the result follows from Theorem \ref{teigF}. }
\end{proof}

\begin{cor}\label{eigN}
Let $k\leq v/2$. Then
\begin{itemize}
\item [\rm(i)] the eigenvalues of $N_{k,k,t}$ are  $$\lambda_j=(-1)^{k-t}{2k-v-1\choose k-j}-\sum_{i=t+1}^k (-1)^{i-t} {k-j\choose i-j}{v-j-i\choose k-i};$$
  \item[\rm(ii)] the eigenvalues of the matrix $N_{k,k,k-1}$ are
$$\lambda_0^{{v \choose 0}},\lambda_1^{{v \choose 1}-{v \choose 0}},\ldots,\lambda_{k-1}^{{v \choose {k-1}}-{v \choose {k-2}}},0^{{v \choose k}-{v \choose {k-1}}},$$
where
$$ \lambda_j=1-{2k-v-1\choose k-j},~~\hbox{for}~j=0,1,\ldots,k-1;$$
  \item[\rm(iii)] with $v=2k$, $\rk N_{k,k,k-1}=\frac{1}{2}{2k\choose k};$
  \item[\rm(iv)]$\rk N_{k,k,k-1}={v \choose {k-1}}$ provided that $k<v/2$.
\end{itemize}
\end{cor}
\begin{proof}
{
\begin{itemize}
\item[\rm(i)] By Corollary~\ref{teigU} and  (\ref{knuth}), the eigenvalues of $N_{k,k,t}=(-1)^t U_{k,k,t,0}$ are  $$\lambda_j=(-1)^{k-t}{2k-v-1\choose k-j}-\sum_{i=t+1}^k (-1)^{i-t} {k-j\choose i-j}{v-j-i\choose k-i}.$$
\item[\rm(ii)] This is an immediate consequence of part (i).
\item[\rm(iii)]
      Setting $v=2k$ in part (ii) yields $\lambda_j=1+{-1 \choose k-j}=1+(-1)^{k-j+1}$ for $j=0,1,\ldots,k-1$. Hence $\lambda_j \neq 0$ if and only if $j+k-1$ is even. Therefore
       $$\rk N_{k,k,k-1}=\sum_{\substack{0\leqslant j\leqslant k-1\\
       2\mid j+k-1}}\left( {2k \choose j}-{2k \choose j-1}\right).$$ The result now follows from the identities
      $\sum_j {2k \choose j}=2^{2k-1}$ and $2\sum_j {2k \choose j-1}=2^{2k}-{2k \choose k}$, where $j$ runs over the same set as in the above sum. (We remark that a direct proof is obtained simply by considering the entries of $N_{k,k,k-1}$.)
  \item[\rm(iv)]It is easily seen that in this case $\lambda_j \neq 0$ for all $j=0,1,\ldots,k-1$, thus $\rk N_{k,k,k-1}=\sum_{j=0}^{k-1}\left( {v \choose j}-{v \choose j-1}\right)={v \choose k-1}$.
\end{itemize}}
\end{proof}

\begin{example}
Considering Remark \ref{N13-14}, we give eigenvalues of the matrices $N_{7,7,6}$ with $v=14$ and $N_{6,6,5}$ with $v=13$. For the first matrix, $\lambda_j=1-{-1\choose 7-j}=1+(-1)^j,$ for $j=0,1,\ldots,6$. Thus the set of eigenvalues is
$\left\{2^{1716},0^{1716}\right\}$ and the rank of this matrix is $1716$.
For second the matrix, $\lambda_j=1-{-2\choose 6-j}=1+(-1)^{j+1}(7-j)$,
for $j=0,1,\ldots,6$. Thus the set of eigenvalues is
$$\left\{(-6)^{1},7^{12},(-4)^{65},5^{208},(-2)^{429},3^{572},0^{429}\right\}$$
and the rank of this matrix is $1287$.
\end{example}

\begin{thm}\label{teigFsk} Let $0\leq t\leq s\leq k\leq v/2$. Consider $F_{s,k,t}(z)$ as a matrix with entries in the
field of rational functions $\mathbb{R}(z)$. Then the eigenvalues
of the matrix $W^\top_{s,k}F_{s,k,t}(z)$ are
$$\alpha_0(z)^{{v \choose 0}},\alpha_1(z)^{{v \choose 1}-{v \choose 0}},\ldots,\alpha_{t}(z)^{{v \choose {t}}-{v \choose {t-1}}},0^{{v \choose k}-{v \choose {t}}},$$
where
   $$\alpha_j(z)= L_{k,s}\mu_j(z) = (-1)^{k+s}\sum_{i=j}^{t}{k-j \choose i-j}{v-j-i\choose k-i}{i-s-1 \choose k-s}z^i.$$
\end{thm}
\begin{proof}
{Again we remark that in Theorem \ref{teigF} the eigenspace of a
given eigenvalue of $F_{k,k,t}(z)$ has a basis independent of
$z$. From this and the equation
$W^\top_{s,k}F_{s,k,t}=L_{k,s}F_{k,k,t}$ (see the proof of
Theorem~\ref{WFUU}) it turns out that the eigenvalues of
$W^\top_{s,k}F_{s,k,t}$ are of the following form:
$$(L_{k,s}\mu_0(z))^{v \choose 0},(L_{k,s}\mu_1(z))^{{v \choose 1}-{v \choose 0}},\ldots,(L_{k,s}\mu_t(z))^{{v \choose t}-{v \choose {t-1}}},0^{{v \choose k}-{v \choose t}},$$
which yields the result. }
\end{proof}

\begin{cor}\label{Urank}
Let $0\leq s\leq k\leq v/2$ and let $\alpha'_j(z)$ be the polynomial obtained from $\alpha_j(z)$ in the previous theorem by setting $t=k$. Then the eigenvalues of the matrix
$W_{s,k}^\top U_{s,k,\ell}$ are
$$\tau_0^{{v \choose 0}},\tau_1^{{v \choose 1}-{v \choose 0}},\ldots,\tau_{s}^{{v \choose s}-{v \choose s-1}},0^{{v \choose k}-{v \choose s}},$$
where
 $$  \tau_j= \frac{D^\ell}{\ell !} \alpha'_j(z)\big|_{z=-1}=(-1)^{k+s+\ell}\sum_{i=\min(j,\ell)}^k (-1)^i {i\choose \ell}{k-j\choose i-j}{v-j-i\choose
k-i}{i-s-1 \choose k-s},$$
for $j=0,1,\ldots,s$. Hence
$$\rk U_{s,k,\ell}=\sum_{\substack{0\leqslant j\leqslant s \\\tau_j\ne0}}  \left( {v \choose j}-{v \choose j-1}\right).$$
\end{cor}

\noindent{\bf Acknowledgments.} The authors are grateful to the referees whose comments greatly improved the presentation
of the paper. The research of the first author was in
part supported by a grant from IPM (No. 89050046). The research of the second author was in
part supported by a grant from IPM (No. 90050117).

\end{document}